\documentstyle[12pt]{amsart}
\begin{document}
  \title{Polyharmonic maps into the Euclidean space}
  \title[Polyharmonic maps into the Euclidean space]
  {Polyharmonic maps into the Euclidean space}
\author{Nobumitsu Nakauchi}
  \address{Graduate School of Science and Engineering, \newline
  Yamaguchi University, 
  Yamaguchi, 753-8512, Japan}
  \email{nakauchi@@yamaguchi-u.ac.jp}
   \author{Hajime Urakawa}
  \address{Division of Mathematics, Graduate School of Information Sciences, Tohoku University, Aoba 6-3-09, Sendai, 980-8579, Japan}
  \curraddr{Institute for International Education, 
  Tohoku University, Kawauchi 41, Sendai 980-8576, Japan}
  \email{urakawa@@math.is.tohoku.ac.jp}
  \thanks{
  Supported by the Grant-in-Aid for the Scientific Reserch, (C) No. 21540207, Japan Society for the Promotion of Science. 
  }
\maketitle
\begin{abstract}
  We study polyharmonic ($k$-harmonic) maps between Riemannian manifolds with finite $j$-energies $(j=1,\cdots,2k-2)$. We show that if the domain is complete and the target is the Euclidean space, then such a map is harmonic. 
    \end{abstract}
\numberwithin{equation}{section}
\theoremstyle{plain}
\newtheorem{df}{Definition}[section]
\newtheorem{th}[df]{Theorem}
\newtheorem{prop}[df]{Proposition}
\newtheorem{lem}[df]{Lemma}
\newtheorem{cor}[df]{Corollary}
\newtheorem{rem}[df]{Remark}
\section{Introduction}
This paper is an extension of our previous work (\cite{NUG}) to polyharmonic maps.  
Harmonic maps play a central role in geometry;\,they are critical points of the energy functional 
$E(\varphi)=\frac12\int_M\vert d\varphi\vert^2\,v_g$ 
for smooth maps $\varphi$ of $(M,g)$ into $(N,h)$. The Euler-Lagrange equations are given by the vanishing of the tension filed 
$\tau(\varphi)$. 
In 1983, J. Eells and L. Lemaire \cite{EL1} extended the notion of harmonic map to  
polyharmonic map, which are, 
by definition, 
critical points of the $k$-energy ($k\geq 2$) 
\begin{equation}
E_k(\varphi)=\frac12\int_M
\vert (d+\delta)^k \varphi\vert^2\,v_g.
\end{equation}
After G.Y. Jiang \cite{J} studied the first and second variation formulas of $E_2$ ($k=2$), 
extensive studies in this area have been done
(for instance, see 
\cite{BFO}, \cite{CMP}, \cite{LO}, \cite{LO2},  \cite{MO1}, \cite{O1},  \cite{S1},
\cite{IIU2}, \cite{IIU},  \cite{II}, 
 etc.). Notice that harmonic maps are always polyharmonic by definition. 
\par
For harmonic maps, it is well known that: 
\par
{\em If a domain manifold $(M,g)$ is complete and has non-negative Ricci curvature, and the sectional curvature of a target manifold $(N,h)$ is non-positive, then 
every energy finite harmonic 
map is a constant map} (cf. \cite{SY}). 
\vskip0.3cm\par
In our previous paper, we showed that 
\begin{th} $($\cite{NUG}$)$
Let 
$(M,g)$ be a complete Riemannian manifold, and 
the curvature of $(N,h)$ is non-positive. Then,   \par
$(1)$ every biharmonic map 
$\varphi:\,(M,g)\rightarrow (N,h)$ with 
finite energy and 
finite bienergy 
must be harmonic. 
\par
$(2)$ In the case ${\rm Vol}(M,g)=\infty$,   
every biharmonic map 
$\varphi:\,(M,g)\rightarrow (N,h)$ with finite bienergy 
is harmonic. 
\end{th}
\vskip0.3cm\par
Now, in this paper, we want to extend it to $k$-harmonic maps 
$(k\geq2$). Indeed,  
we will show 
\begin{th} (Theorems 2.4 and 3.1)\quad
Let $(M,g)$ be a complete 
Riemannian manifold, and 
$(N,h)$, the $n$-dimensional Euclidean space. 
Then, 
\par 
$(1)$ 
every 
$k$-harmonic map $\varphi:\,(M,g)\rightarrow (N,h)$ $(k\geq2)$ with 
finite  $j$-energies for all $j=1,2,\, \cdots,2k-2$, must be harmonic. 
\par
$(2)$ In the case of \,\,$\mbox{\rm Vol}(M,g)=\infty$, 
every $k$-harmonic map $\varphi:\,(M,g)\rightarrow (N,h)$ 
with finite $j$-energy for all $j=2,4,\,\cdots,2k-2$, is harmonic. 
\end{th}
\vskip0.3cm\par
Theorem 1.2 
gives an affirmative answer to the generalized B.Y. Chen's conjecture (cf. \cite{CMP})
on $k$-harmonic maps $(k\geq 2$) 
under the $L^2$-conditions. 
\vskip0.3cm\par
{\bf Acknowledgement.}\quad We express our gratitude to Dr. Shun Maeta who gave
 valuable comments in the first draft. 
\vskip0.6cm\par
\section{Preliminaries and statement of main theorem}
In this section, we prepare materials for the first variational formula for the biharmonic maps. 
Let us recall the definition of a harmonic map $\varphi:\,(M,g)\rightarrow (N,h)$, of a compact Riemannian manifold $(M,g)$ into another Riemannian manifold $(N,h)$, 
which is an extremal 
of the {\em energy functional} defined by 
$$
E(\varphi)=\int_Me(\varphi)\,v_g, 
$$
where $e(\varphi):=\frac12\vert d\varphi\vert^2$ is called the energy density 
of $\varphi$.  
That is, for any variation $\{\varphi_t\}$ of $\varphi$ with 
$\varphi_0=\varphi$, 
\begin{equation}
\frac{d}{dt}\bigg\vert_{t=0}E(\varphi_t)=-\int_Mh(\tau(\varphi),V)v_g=0,
\end{equation}
where $V\in \Gamma(\varphi^{-1}TN)$ is a variation vector field along $\varphi$ which is given by 
$V(x)=\frac{d}{dt}\vert_{t=0}\varphi_t(x)\in T_{\varphi(x)}N$, 
$(x\in M)$, 
and  the {\em tension field} is given by 
$\tau(\varphi)
=\sum_{i=1}^mB(\varphi)(e_i,e_i)\in \Gamma(\varphi^{-1}TN)$, 
where 
$\{e_i\}_{i=1}^m$ is a locally defined frame field on $(M,g)$, 
and $B(\varphi)$ is the second fundamental form of $\varphi$ 
defined by 
\begin{align}
B(\varphi)(X,Y)&=(\widetilde{\nabla}d\varphi)(X,Y)\nonumber\\
&=(\widetilde{\nabla}_Xd\varphi)(Y)\nonumber\\
&=\overline{\nabla}_X(d\varphi(Y))-d\varphi(\nabla_XY),
\end{align}
for all vector fields $X, Y\in {\frak X}(M)$. 
Here, 
$\nabla$, and
$\nabla^N$, 
 are the Levi-Civita connections of $(M,g)$, $(N,h)$, respectively, and 
$\overline{\nabla}$, and $\widetilde{\nabla}$ are the induced ones on $\varphi^{-1}TN$, and $T^{\ast}M\otimes \varphi^{-1}TN$, respectively. By (2.1), $\varphi$ is harmonic if and only if $\tau(\varphi)=0$. 
\par
The second variation formula is given as follows. Assume that 
$\varphi$ is harmonic. 
Then, 
\begin{equation}
\frac{d^2}{dt^2}\bigg\vert_{t=0}E(\varphi_t)
=\int_Mh(J(V),V)v_g, 
\end{equation}
where 
$J$ is an elliptic differential operator, called the
{\em Jacobi operator}  acting on 
$\Gamma(\varphi^{-1}TN)$ given by 
\begin{equation}
J(V)=\overline{\Delta}V-{\mathcal R}(V),
\end{equation}
where 
$\overline{\Delta}V=\overline{\nabla}^{\ast}\overline{\nabla}V
=-\sum_{i=1}^m\{
\overline{\nabla}_{e_i}\overline{\nabla}_{e_i}V-\overline{\nabla}_{\nabla_{e_i}e_i}V
\}$ 
is the {\em rough Laplacian} and 
${\mathcal R}$ is a linear operator on $\Gamma(\varphi^{-1}TN)$
given by 
${\mathcal R}(V)=
\sum_{i=1}^mR^N(V,d\varphi(e_i))d\varphi(e_i)$,
and $R^N$ is the curvature tensor of $(N,h)$ given by 
$R^N(U,V)=\nabla^N{}_U\nabla^N{}_V-\nabla^N{}_V\nabla^N{}_U-\nabla^N{}_{[U,V]}$ for $U,\,V\in {\frak X}(N)$.   
\par
J. Eells and L. Lemaire \cite{EL1} proposed polyharmonic ($k$-harmonic) maps and 
Jiang \cite{J} studied the first and second variation formulas for biharmonic maps. Let us consider the {\em bienergy functional} 
defined by 
\begin{equation}
E_2(\varphi)=\frac12\int_M\vert\tau(\varphi)\vert ^2v_g, 
\end{equation}
where 
$\vert V\vert^2=h(V,V)$, $V\in \Gamma(\varphi^{-1}TN)$.  
The first variation formula of the bienergy functional 
is given by
\begin{equation}
\frac{d}{dt}\bigg\vert_{t=0}E_2(\varphi_t)
=-\int_Mh(\tau_2(\varphi),V)v_g.
\end{equation}
Here, 
\begin{equation}
\tau_2(\varphi)
:=J(\tau(\varphi))=\overline{\Delta}(\tau(\varphi))-{\mathcal R}(\tau(\varphi)),
\end{equation}
which is called the {\em bitension field} of $\varphi$, and 
$J$ is given in $(2.4)$.  
\par
A smooth map $\varphi$ of $(M,g)$ into $(N,h)$ is said to be 
{\em biharmonic} if 
$\tau_2(\varphi)=0$. 
\vskip0.3cm\par
Now let us recall the definition of the $k$-energy $E_k(\varphi)$ ($k\geq 2$):
\begin{df}  
The {\em $k$-energy $E_k(\varphi)$} ($k\geq 2$) is defined formally $($\cite{EL2}$)$ by 
\begin{align}
E_k(\varphi):=\frac12\int_M\vert (d+\delta)^k\varphi\vert^2\,v_g
\end{align}
for every smooth map $\varphi\in C^{\infty}(M,N)$. 
Then, it is given $($\cite{IIU2}, p. 270, Lemma 40$)$ by the following formula:
\begin{equation}
E_k(\varphi)=\left\{
\begin{aligned}
&\frac12
\int_M\vert W^{\ell}_{\varphi}\vert^2\,v_g\qquad (\mbox{\rm if}\,\,k\,\,\mbox{\rm is even, say}\,\,2\,\ell),\\
&\frac12\int_M\vert \overline{\nabla}\,W^{\ell}_{\varphi}\vert^2\,v_g
\qquad (\mbox{\rm if}\,\,k\,\,\mbox{\rm is odd, say}\,\,2\ell+1). 
\end{aligned}
\right.
\end{equation}
Here, $W^{\ell}_{\varphi}$ is given as, by definition,  
\begin{equation}
W^{\ell}_{\varphi}:=\underbrace{\overline{\Delta}\cdots\overline{\Delta}}_{\ell-1}
\,\tau(\varphi)\in \Gamma(\varphi^{-1}TN).
\end{equation}
For $k=1$, that is, $\ell=0$, we define $W^0_{\varphi}=\varphi$, also. 
\end{df}
\vskip0.6cm\par
Then, the definition and the first variation formula for the $k$-energy $E_k$ are given 
as follows: 
\begin{df}\,\, $(${\bf $k$-harmonic map}$)$ 
For each $k=2,3,\cdots$, and a smooth map 
$\varphi:\,(M,g)\rightarrow (N,h)$, is $k$-{\em harmonic} if 
\begin{align}
\frac{d}{dt}\bigg\vert_{t=0}E_k(\varphi_t)=0
\end{align}
for every smooth variation $\varphi_t:\,M\rightarrow N$ ($-\epsilon<t<\epsilon$) 
with $\varphi_0=\varphi$. 
\end{df}
\vskip0.6cm\par
Then, we have (\cite{IIU2}, p.269, Theorem 39) 
\begin{th} $(${\bf The first variation formula of the $k$-energy}$)$ \quad 
Assume that $(N,h)=({\mathbb R}^n,h_{{\mathbb R}^n})$ is the $n$-dimensional Euclidean space. 
For every $k=2,3,\cdots$, it holds that 
\begin{align}
\frac{d}{dt}\bigg\vert_{t=0}E_k(\varphi_t)=-\int_M
\langle \tau_k(\varphi),V\rangle\,v_g,
\end{align}
where $V$ is a variation vector field given by 
$V(x)=\frac{d}{dt}\big\vert_{t=0}\varphi_t(x)\in T_{\varphi(x)}N$ ($x\in M$).
The {\em $k$-tension field} $\tau_k(\varphi)$ is given by 
\begin{align}
\tau_k(\varphi)=J(W_{\varphi}^{k-1})
=\overline{\Delta}(W^{k-1}_{\varphi}),
\end{align}
where $W^{k-1}_{\varphi}=\underbrace{\overline{\Delta}\cdots\overline{\Delta}}_{k-2}\,\tau(\varphi)\in \Gamma(\varphi^{-1}TN)$. 
\par
Thus, $\varphi:\,(M,g)\rightarrow (N,h)$ is $k$-harmonic if and only if 
$\overline{\Delta}^{k-1}\tau(\varphi)=0$ which is equivalent to $W_{\varphi}^k=0$. 
\end{th}
\vskip0.6cm\par
The formula (143) of the $k$-tension field $\tau_k(\varphi)$ in Theorem 39 (p.269, \cite{IIU2}) is true only for the case that the target space $(N,h)$ is the $n$-dimensional Euclidean space $(N,h)=({\mathbb R}^n,h_{{\mathbb R}^n})$. 
\par
Here, we denote by 
$\overline{\nabla}W^{\ell}_{\varphi}=\overline{\nabla}\varphi=d\varphi$ 
for $\ell=0$, and $k=2\ell+1=1$,
$$
E_1(\varphi)=\frac12\int_M\vert d\varphi\vert^2\,v_g. 
$$ 
\vskip0.6cm\par
Then, we can state our main theorem. 
\begin{th} 
$(${\bf Main theorem}$)$ 
Assume that the domain manifold $(M,g)$ is a complete Riemannian manifold, and the target space $(N,h)$ is the $n$-dimensional Euclidean space. 
Let $\varphi:\,(M,g)\rightarrow (N,h)$ be a $k$-harmonic map ($k\geq 2$). 
 Assume that 
 \begin{align}
 &(1) \,\,\mbox{$E_j(\varphi)<\infty$ for all}\,\, j= 2,4,\cdots,2k-2,\,\,\mbox{and}\nonumber\\
 &(2) \,\, \mbox{either} \nonumber\\
 &\qquad\quad \mbox{$E_j(\varphi)<\infty$ for all}\,\, j=1,3,\cdots, 2k-3,\,\,\mbox{or}\nonumber\\
 &\qquad\qquad 
 \mbox{\rm Vol}(M,g)=\infty. \nonumber
 \end{align}
Then, $\varphi:\,(M,g)\rightarrow (N,h)$  
is harmonic. 
\end{th}
In the case of the $n$-dimensional Euclidean space 
$(N,h)=({\mathbb R}^n,h_{{\mathbb R}^n})$, 
Theorem 2.4 and the following Theorem 3.1 are natural extensions of our previous theorem in \cite{NUG} which is: 
\begin{th} 
Assume that $(M,g)$ is complete and the sectional curvature of $(N,h)$ 
is non-positive. \par
$(1)$ Every biharmonic map 
$\varphi:\,(M,g)\rightarrow (N,h)$ with 
finite energy $E(\varphi)<\infty$ and 
finite bienergy $E_2(\varphi)<\infty$, 
is harmonic. \par
$(2)$ In the case ${\rm Vol}(M,g)=\infty$,  
every biharmonic map 
$\varphi:\,(M,g)\rightarrow (N,h)$ with finite bienergy $E_2(\varphi)<\infty$,
is harmonic. 
\end{th}
\vskip0.6cm\par
\section{The iteration proposition.}
By virtue of (2.9), we have to notice the the energy conditions in (1) and (2) of Theorem 2.4: 
\par
Indeed, the condition which $E_j(\varphi)<\infty$ for all $j=2,4,\cdots, 2k-2$ in (1) of Theorem 2.4 
is equivalent to that 
\begin{align}
\int_M\vert W^j_{\varphi}\vert^2\,v_g<\infty\qquad (j=1,2,\cdots,k-1), 
\end{align}
and the condition which $E_j(\varphi)<\infty$ for all $j=1,3,\cdots,2k-3$ in (2) of Theorem 2.4 
is equivalent to that 
\begin{align}
\int_M \vert\overline{\nabla} W^j_{\varphi}\vert^2\,v_g<\infty 
\qquad (j=0,1,\cdots,k-2).  
\end{align}
Therefore, to show Theorem 2.4, 
we only have to prove the following theorem:
\begin{th}
Assume that the domain manifold $(M,g)$ is a complete Riemannian manifold, and the target space $(N,h)$ is the $n$-dimensional Euclidean space. 
Let $\varphi:\,(M,g)\rightarrow (N,h)$ be a $k$-harmonic map. 
\par
Assume that  
$$
(1) \qquad \int_M\vert W^j_{\varphi}\vert^2\,v_g<\infty\,\,\mbox{for all $j=1,2,\cdots,k-1$}, \,\,\mbox{and}\qquad\qquad\qquad
$$
$(2)$\quad \mbox{either} \,\,
 $$\int_M\vert\overline{\nabla}\,W^j_{\varphi}\vert^2\,v_g<\infty\,\,
\mbox{for all} \,\,j=0,1,\cdots,k-2,\,\,\mbox{or}$$  \,\,
$$
\mbox{\rm Vol}(M,g)=\infty. \qquad\qquad\qquad\qquad\qquad\qquad\qquad
$$
\par
Then, $\varphi:\,(M,g)\rightarrow (N,h)$ is harmonic. 
\end{th}
To prove Theorem 3.1 whose proof will be given in the next section, we need the following iteration proposition: 
\begin{prop} $($the iteration method$)$
Let $(M,g)$ be a complete Rienannian manifold, and $(N,h)$, an arbitrary Riemannian manifold. 
 Let $\varphi:\,(M,g)\rightarrow (N,h)$ be an arbitrary $C^{\infty}$ map 
 satisfying that 
 for some $j\geq 2$, 
 \par
 \begin{equation}
 W_{\varphi}^j=0. 
 \end{equation}
 If we assume the following two conditions: 
 \begin{equation}
 \left\{
 \begin{aligned}
 \mbox{$(1)$} \quad &\int_M\vert W_{\varphi}^{j-1}\vert^2\,v_g<\infty, \mbox{\it and}\\ 
 \mbox{$(2)$} \quad &\mbox{either $\int_M\vert\overline{\nabla}\,W^{j-2}_{\varphi}\vert^2\,v_g<\infty$ or 
 $\mbox{\rm Vol}(M,g)=\infty$,} 
 \end{aligned}
 \right.
 \end{equation}
 \par\noindent
 then, we have 
 \begin{equation}
 W^{j-1}_{\varphi}=0. 
 \end{equation}
\end{prop}
\begin{rem}
Under the assumptions (3.2), if we have $W^k_{\varphi}=0$ for some $k\geq 2$, 
then we have automatically, $W^1_{\varphi}=\tau(\varphi)=0$, i.e., $\varphi$ is harmonic. 
\end{rem}
\vskip0.6cm\par
In this section, we give a proof of 
Proposition 3.2 which consists of 
four steps. 
\par
({\it The first step}) \quad 
For a fixed point $x_0\in M$, and 
for every 
$0<r<\infty$, 
we first take a cut-off  $C^{\infty}$ function $\eta$ on $M$ 
(for instance, see \cite{K}) satisfying that 
\begin{equation}
\left\{
\begin{aligned}
0\leq &\eta(x)\leq 1\quad (x\in M),\\
\eta(x)&=1\qquad\quad (x\in B_r(x_0)),\\
\eta(x)&=0\qquad\quad (x\not\in B_{2r}(x_0)),\\
\vert\nabla\eta\vert&\leq\frac{2}{r}
\qquad\,\,\, (x\in M).
\end{aligned}
\right.
\end{equation}
\par
\vskip0.6cm\par
({\it The second step}) \quad 
Notice that (3.3) is equivalent to that 
\begin{equation}
\overline{\Delta} \,W^{j-1}_{\varphi}=0
\end{equation}
because of $W^j_{\varphi}=\overline{\Delta}\,W^{j-1}_{\varphi}$.  
\par
Then, we have 
\begin{align}
0&=\int_M\langle\eta^2\,W^{j-1}_{\varphi},\overline{\Delta}\,W^{j-1}_{\varphi}\rangle\,v_g\nonumber\\
&=\int_M\sum_{i=1}^m\langle \overline{\nabla}_{e_i}(\eta^2\,W^{j-1}_{\varphi}),
\overline{\nabla}_{e_i}W^{j-1}_{\varphi}\rangle\,v_g\nonumber\\
&=\int_M
\eta^2\sum_{i=1}^m\vert\overline{\nabla}_{e_i}W^{j-1}_{\varphi}\vert^2v_g+
2\int_M\sum_{i=1}^m\eta\,e_i(\eta)\langle W^{j-1}_{\varphi},\overline{\nabla}_{e_i}W^{j-1}_{\varphi}\rangle\,v_g.
\end{align}
By moving the second term in the last equality of (3.8) to the left hand side,  
we have 
\begin{align}
\int_M\eta^2\,\sum_{i=1}^m 
\vert\overline{\nabla}_{e_i} W^{j-1}_{\varphi}\vert^2
&=
-2\int_M\sum_{i=1}^m\langle \eta\,\overline{\nabla}_{e_i} W^{j-1}_{\varphi},
e_i(\eta)\,W^{j-1}_{\varphi}\rangle\,v_g\nonumber\\
&=-2\int_M
\sum_{i=1}^m\langle S_i,T_i\rangle\,v_g, 
\end{align}
where we put 
$S_i:=\eta\,\overline{\nabla}_{e_i}W^{j-1}_{\varphi}$, and 
$T_i:=e_i(\eta)\,W^{j-1}_{\varphi}$ ($i=1\,\cdots,m$).  
 \par
 Now let recall the following 
 inequality: 
 \begin{equation}
 \pm 2\,\langle S_i,T_i\rangle
 \leq \epsilon \vert S_i\vert ^2+\frac{1}{\epsilon}\vert T_i\vert^2
 \end{equation} 
 for all positive $\epsilon>0$ because of the inequality 
 $
 0\leq \vert \sqrt{\epsilon}\,S_i\pm\frac{1}{\sqrt{\epsilon}}\,T_i\vert^2.  
 $
 Therefore, for (3.10), we obtain 
\begin{align}
-2\int_M\sum_{i=1}^m\langle S_i,T_i\rangle\,v_g
\leq
\epsilon\int_M
\sum_{i=1}^m\vert S_i\vert^2\,v_g
+\frac{1}{\epsilon}\int_M\sum_{i=1}^m\vert T_i\vert^2\,v_g. 
\end{align} 
If we put $\epsilon=\frac12$, we obtain, 
by (3.9) and (3.11), 
\begin{align}
\int_M
\eta^2\sum_{i=1}^m\vert\overline{\nabla}_{e_i}W^{j-1}_{\varphi}\vert^2\,v_g
&\leq \frac12 
\int_M\sum_{i=1}^m\eta^2\,\vert\overline{\nabla}_{e_i}W^{j-1}_{\varphi}\vert^2\,v_g
\nonumber\\
&\qquad+2\int_M\sum_{i=1}^me_i(\eta)^2\,\vert W^{j-1}_{\varphi}\vert^2\,v_g.
\end{align}
Thus, by (3.12) and (3.6), we obtain 
\begin{align}
\int_M\eta^2\sum_{i=1}^m\vert\overline{\nabla}_{e_i}W^{j-1}_{\varphi}\vert^2\,v_g
&\leq 
4\int_M\vert\nabla\eta\vert^2\,\vert W^{j-1}_{\varphi}\vert^2\,v_g\nonumber\\
&\leq\frac{16}{r^2}\int_M\vert W^{j-1}_{\varphi}\vert^2\,v_g.
\end{align}
\par
({\it The third step}) \quad By definition of $\eta$ in the first step, (3.13) turns out that 
\begin{align}
\int_{B_r(x_0)}\vert\overline{\nabla}\,W^{j-1}_{\varphi}\vert^2\,v_g\leq 
\frac{16}{r^2}\int_M\vert W^{j-1}_{\varphi}\vert^2\,v_g.
\end{align}
Here, recall our assumption that $(M,g)$ is complete and non-compact, and 
(1) $\int_M\vert W^{j-1}_{\varphi}\vert^2\,v_g<\infty$. 
When we tend $r\rightarrow \infty$, the right hand side of (3.12) goes to zero, and the left hand side of (3.12) 
goes to $\int_M\vert \overline{\nabla}W^{j-1}_{\varphi}\vert^2\,v_g$. 
Thus, we obtain 
$$
0\leq \int_M\vert \overline{\nabla}W^{j-1}_{\varphi}\vert^2\,v_g\leq 0,
$$
which implies that 
\begin{align}
\overline{\nabla}\,W^{j-1}_{\varphi}=0
\end{align}
everywhere on $M$. 
\par
({\it The fourth step}) \quad 
(a) In the case that $\int_M\vert\overline{\nabla}\,W^{j-2}_{\varphi}\vert^2\,v_g<\infty$, let us define 
a smooth $1$-form $\alpha$ on $M$ by 
\begin{equation}
\alpha(X):=\langle W^{j-1}_{\varphi},\overline{\nabla}_XW^{j-2}_{\varphi}\rangle\qquad (X\in {\frak X}(M).
\end{equation}
Then, we have: 
\begin{align}
\mbox{\rm div}(\alpha)=-\vert W^{j-1}_{\varphi}\vert^2. 
\end{align}
Because we have 
\begin{align}
\mbox{\rm div}(\alpha)&=\sum_{i=1}^m(\nabla_{e_i}\alpha)(e_i)\nonumber\\
&=\sum_{i=1}^m
\{e_i(\alpha(e_i))-\alpha(\nabla_{e_i}e_i)
\}\nonumber\\
&=\sum_{i=1}^m\bigg\{
e_i\,\big(\langle\,W^{j-1}_{\varphi},\overline{\nabla}_{e_i}W^{j-2}_{\varphi}\rangle \big)
-\langle W^{j-1}_{\varphi},\overline{\nabla}_{\nabla_{e_i}e_i}W^{j-2}_{\varphi}
\rangle
\bigg\}\nonumber\\
&=\sum_{i=1}^m
\bigg\{
\langle \overline{\nabla}_{e_i}W^{j-1}_{\varphi},\overline{\nabla}_{e_i}W^{j-2}_{\varphi} \rangle
+\langle W^{j-1}_{\varphi},\overline{\nabla}_{e_i}\overline{\nabla}_{e_i}W^{j-2}_{\varphi}\rangle 
\nonumber\\
&\qquad\qquad\qquad\qquad
-\langle W^{j-1}_{\varphi},\overline{\nabla}_{\nabla_{e_i}e_i}W^{j-2}_{\varphi}
\bigg\}\nonumber\\
&=\langle W^{j-1}_{\varphi},-\overline{\Delta}W^{j-2}_{\varphi}\rangle 
\qquad(\mbox{because of (3.15) and definition of $\overline{\Delta}$})\nonumber\\
&=-\vert W^{j-1}_{\varphi}\vert^2,
\end{align}
which is (3.17). 
\par
Furthermore, we have 
\begin{align}
\int_M\vert\alpha\vert\,v_g<\infty.
\end{align}
Because we have, by definition of $\alpha$ in (3.16),
\begin{align}
\int_M\vert\alpha\vert\,v_g&=
\int_M\vert 
\langle W^{j-1}_{\varphi},\overline{\nabla}W^{j-2}_{\varphi}\rangle\vert\,v_g\nonumber\\
&\leq \bigg(
\int_M\vert W^{j-1}_{\varphi}\vert^2\,v_g
\bigg)^{\frac12}\bigg(
\int_M\vert\overline{\nabla}W^{j-2}_{\varphi}\vert^2\,v_g
\bigg)^{\frac12}\nonumber\\
&<\infty
\end{align}
because of our assumptions 
$\int_M\vert W^{j-1}_{\varphi}\vert^2\,v_g<\infty$ and 
$\int_M\vert\overline{\nabla}W^{j-2}_{\varphi}\vert^2\,v_g<\infty$. 
Thus, we can apply Gaffney's theorem to this $\alpha$ (cf. \cite{G}, 
and Theorem 4.1 in Appendix in \cite{NUG}). 
 We obtain 
 \begin{align}
 0=\int_M\mbox{\rm div}(\alpha)\,v_g=-\int_M\vert W^{j-1}_{\varphi}\vert^2\,v_g,
 \end{align}
 which implies that 
 $W^{j-1}_{\varphi}=0$. 
 \par
 (b) In the case that $\mbox{\rm Vol}(M,g)=\infty$, 
 we first notice that 
 $\vert W^{j-1}_{\varphi}\vert^2$ is constant on $M$, say $C_0$. 
 Because for every $X\in {\frak X}(M)$, we have 
 \begin{align}
 X\,\vert W^{j-1}_{\varphi}\vert^2=2\,\langle 
 \overline{\nabla}_XW^{j-1}_{\varphi},W^{j-1}_{\varphi}\rangle=0
 \end{align}
 due to (3.15). Then, due to the assumption (1) of Proposition 3.2,  and the above, 
 we obtain 
 \begin{align}
 \infty>\int_M\vert W^{j-1}_{\varphi}\vert^2\,v_g=C_0\,\int_Mv_g=C_0\,\mbox{\rm Vol}(M,g). 
 \end{align} 
 By our assumption that $\mbox{\rm Vol}(M,g)=\infty$, 
 (3.23) implies that $C_0=0$. We obtain 
 $W^{j-1}_{\varphi}\equiv 0$.  
 We obtain Proposition 3.2.
 \qed
 \vskip0.6cm\par
{\it Proof of Theorem 3.1.}
\quad 
We apply Proposition 3.2 to our map $\varphi:\,(M,g)\rightarrow (N,h)$, then the iteration procedure works well 
since $\varphi$ is $k$-harmonic, i.e., $W^{k}_{\varphi}=0$. Then, we have 
$W^{k-1}_{\varphi}=0$, and then we have $W^{k-2}_{\varphi}=0$, etc.  Finally, 
we obtain $\tau(\varphi)=W^1_{\varphi}=0$. 
Thus, $\varphi:\,(M,g)\rightarrow (N,h)$ is harmonic. We obtain Theorem 3.1.
\qed
\vskip1cm\par


\begin{thebibliography}{99}
\bibitem{BE} 
P. Baird and J. Eells, 
\textit{A conservation law for harmonic maps}, 
Lecture Notes in Math., Springer, {\bf 894} (1981), 
1--25. 
\bibitem{BFO} 
P. Baird, A. Fardoun and S. Ouakkas, 
\textit{Liouville-type theorems for biharmonic maps between Riemannian manifolds}, 
Adv. Calc. Var., {\bf 3} (2010), 49--68. 
\bibitem{BW} 
P. Baird and J. Wood, 
\textit{Harmonic Morphisms Between Riemannian Manifolds}, Oxford Science Publication, 2003, Oxford. 
\bibitem{CMP} R. Caddeo, S. Montaldo, P. Piu, 
\textit{On biharmonic maps}, Contemp. Math., \textbf{288} (2001), 286--290. 
\bibitem{C} B.Y. Chen, \textit{Some open problems and conjectures on submanifolds of finite type}, 
Soochow J. Math., \textbf{17} (1991), 169--188. 
\bibitem{EL1} J. Eells, L. Lemaire, 
\textit{A report on harmonic maps}, Bull. London Math. Soc., \textbf{10} (1978), 1--68. 
\bibitem{EL2} J. Eells, L. Lemaire, 
\textit{Selected topics in harmonic maps}, CBMS, 
\textbf{50}, Amer. Math. Soc, 1983. 
\bibitem{EL3} J. Eells, L. Lemaire, 
{\em Another Report on Harmonic Maps}, Bull. London Math. Soc., \textbf{20} (1988), 385--524. 
\bibitem{ES} 
J. Eells and J.H. Sampson, 
{\em Harmonic mappings of Riemannian manifolds}, 
Amer. J. Math., {\bf 86} (1964), 109--160. 
\bibitem{G}
M.P. Gaffney 
{\em A special Stokes' theorem for complete Riemannian manifold}, 
Ann. Math., {\bf 60} (1954), 140--145. 
\bibitem{Gu} 
S. Gudmundsson, 
{\it The Bibliography of Harmonic Morphisms}, {\tt http://matematik.lu.se/\\matematiklu/personal/sigma/harmonic/bibliography.html}
\bibitem{IIU2} 
T. Ichiyama, J. Inoguchi, H. Urakawa, \textit{Biharmonic maps and bi-Yang-Mills fields}, 
Note di Matematica, {\bf 28}, (2009), 233--275. 
\bibitem{IIU} T. Ichiyama, J. Inoguchi, H. Urakawa, \textit{Classifications and isolation phenomena  of biharmonic maps and bi-Yang-Mills fields}, 
Note di Matematica, {\bf 30}, (2010), 15--48. 
\bibitem{II} S. Ishihara, S. Ishikawa, \textit{Notes on relatively harmonic immersions}, Hokkaido Math. J., \textbf{4}(1975), 234--246. 
\bibitem{J} G.Y. Jiang, \textit{2-harmonic maps and their first and second variational formula}, Chinese Ann. Math., \textbf{7A} (1986), 
388--402;  Note di Matematica, {\bf 28} (2009), 209--232.
\bibitem{K} A. Kasue, {\em Riemannian Geometry}, in Japanese, 
Baihu-kan, Tokyo, 2001. 
\bibitem{L} T. Lamm, \textit{Biharmonic map heat flow into manifolds of nonpositive curvature}, 
Calc. Var., \textbf{22} (2005), 421--445. 
\bibitem{LO} E. Loubeau, C. Oniciuc, \textit{The index of biharmonic maps in spheres}, 
Compositio Math., \textbf{141} (2005), 729--745. 
\bibitem{LO2} E. Loubeau and C. Oniciuc, \textit{On the biharmonic and harmonic indices of the Hopf map}, Trans. Amer. Math. Soc., {\bf 359} (2007), 5239--5256. 
\bibitem{LOu}
E. Loubeau and Y-L. Ou, 
\textit{
Biharmonic maps and morphisms from conformal mappings}, 
Tohoku Math. J., {\bf 62} (2010), 55--73. 
\bibitem{MO1} S. Montaldo, C. Oniciuc, \textit{A short survey on biharmonic maps between Riemannian manifolds}, Rev. Un. Mat. Argentina \textbf{47} (2006), 1--22. 
\bibitem{NU1} 
N. Nakauchi and H. Urakawa, 
\textit{Biharmonic hypersurfaces in a Riemannian manifold with non-positive Ricci curvature}, 
Ann. Global Anal. Geom., {\bf 40} (2011), 125--131. 
\bibitem{NU2} 
N. Nakauchi and H. Urakawa, 
\textit{Biharmonic submanifolds in a Riemannian 
manifold with non-positive curvature}, Results in Math., \textbf{63} (2013), 467--474. 
\bibitem{NUG} 
N. Nakauchi, H. Urakawa and 
S. Gudmundsson, \textit{Biharmonic maps into a Riemannian manifold of non-positive curvature}, 
to appear in Geometriae Dedicata, 2013. 
\bibitem{O1} C. Oniciuc, \textit{On the second variation formula for biharmonic maps to a sphere}, 
Publ. Math. Debrecen., \textbf{67} (2005), 285--303. 
\bibitem{OT} 
Ye-Lin Ou and Liang Tang, 
{\em The generalized Chen's conjecture on biharmonic submanifolds is false}, 
arXiv: 1006.1838v1.  
\bibitem{S1} T. Sasahara, \textit{Legendre surfaces in Sasakian space forms whose mean curvature vectors are eigenvectors}, Publ. Math. 
Debrecen, \textbf{67} (2005), 285--303. 
\bibitem{SY} R. Schoen and S.T. Yau, 
\textit{Harmonic maps and the topology of stable hypersurfaces and manifolds with non-negative Ricci curvature}, Comment. Math. Helv. {\bf 51} (1976), 333--341. 
\bibitem{Y} 
S.T. Yau, 
\textit{Some function-theoretic properties of complete Riemannian manifold and their applications to geometry}, 
Indiana Univ. Math. J., {\bf 25} (1976), 659--670. 
\bibitem{WO} 
Z-P Wang and Y-L Ou, Biharmonic Riemannian submersions from 3-manifolds, Math. Z., {\bf 269} (2011), 917--925.  
\end{thebibliography}
\end{document}